\documentclass{amsart}
\usepackage{color,graphicx,wrapfig}
\title{The Parabolic Two-Phase Membrane Problem: Regularity in Higher Dimensions}
\author[H. Shahgholian ]{Henrik Shahgholian}
\address{Department of Mathematics, Royal Institute of Technology,
100~44  Stockholm, Sweden}
\email{henriksh@math.kth.se}
\author[N. Uraltseva ]{Nina Uraltseva}
\address{St. Petersburg State University,
Department of Mathematics and Mechanics,
198504, St. Petersburg, Staryi Petergof,
Universitetsky Pr. 28}
\email{uraltsev@pdmi.ras.ru}
\author[G.S. Weiss]{Georg S. Weiss}
\address{Graduate School of Mathematical Sciences,
University of Tokyo, 3-8-1 Komaba, Meguro-ku, Tokyo-to, 153-8914 Japan}
\email{gw@ms.u-tokyo.ac.jp}
\thanks{$2000$ {\it Mathematics Subject Classification.\/} Primary
35R35, Secondary 35J60.}
\thanks{{\it Key words and phrases.\/} Free boundary,
singular point, branch point, membrane, obstacle problem, regularity, global solution, blow-up.}
\thanks{H. Shahgholian has been partially supported  by the Swedish Research Council.
N. Uraltseva has been partially 
supported by Russian Foundation of Basic Research
(grant number 05-01-01063).
G.S. Weiss has been partially supported by a Grant-in-Aid for Scientific Research,
Ministry of Education, Japan}
\thanks{G.S. Weiss and N. Uraltseva wish to thank Swedish STINT and the Royal Swedish Academy of Science, respectively,
 for visiting appointments to the Royal Inst. of Technology in Stockholm. }
\date{}

\theoremstyle{plain}
\newtheorem{theorem}{Theorem}[section]
\newtheorem{lemma}[theorem]{Lemma}
\newtheorem{proposition}[theorem]{Proposition}
\newtheorem{corollary}[theorem]{Corollary} \theoremstyle{definition}
 \theoremstyle{example}

\theoremstyle{definition}

\numberwithin{equation}{section}

\def\R{{\bf R}}

\def\C{{\bf C}}

\def\H{{\bf H}}

\def\dist{\hbox{\rm dist}}
\def\pardist{\hbox{\rm pardist}}

\begin{document}
\maketitle
\begin{abstract}
For the parabolic obstacle-problem-like equation
$$\Delta u - \partial_t u = \lambda_+ \chi_{\{u>0\}} \> - \>
\lambda_- \chi_{\{u<0\}}\> ,$$ where $\lambda_+$ and $\lambda_-$ are
positive Lipschitz functions, we prove in arbitrary finite dimension that the
free boundary $\partial\{ u>0\} \cup\partial\{ u<0\}$ is in a neighborhood of each ``branch point'' the
union of two Lipschitz graphs that are continuously differentiable
with respect to the space variables.
The result extends the elliptic paper \cite{imrn} to the parabolic case. 
The result is optimal in the sense
that the graphs are in general not better than Lipschitz, as shown by a counter-example.
\end{abstract}
\tableofcontents
\section{Introduction}
\subsection{Background and main result}
\begin{figure}
\begin{center}
\input{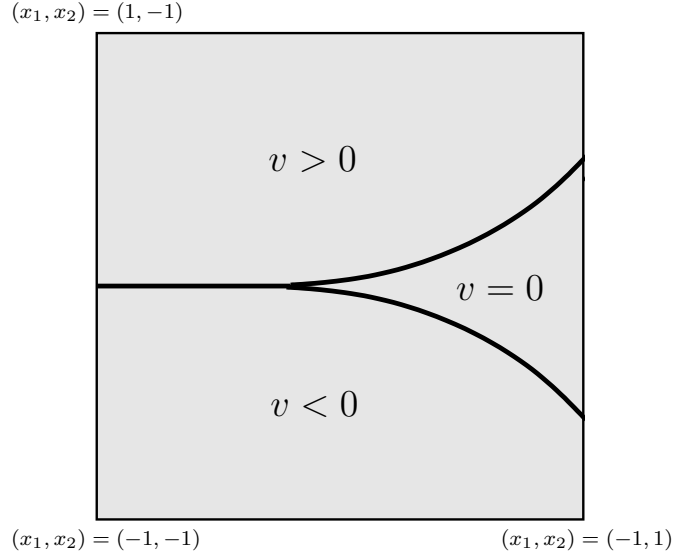}
\end{center}
\caption{Example of a Stationary Branch Point}\label{branchfig}
\end{figure}

In this paper we study the regularity of the
parabolic obstacle-problem-like equation
\begin{equation}\label{obst}
\Delta u - \partial_t u = \lambda_+ \chi_{\{u>0\}} \> - \> \lambda_-
\chi_{\{u<0\}}\;    \qquad \hbox{in } (0,T)\times \Omega ,
\end{equation}
where $T<+\infty, \lambda_+> 0, \lambda_->0$ are Lipschitz functions and $\Omega \subset \R^n$ is a
given domain. The problem arises as limiting case in the model
of temperature control through the interior described
in \cite[2.3.2]{duvaut} as $h_1,h_2\to 0$.\\
We are interested in the regularity of the free boundary $\partial\{ u>0\} \cup\partial\{ u<0\}$.
As the one-phase case (i.e. the case of a non-negative or non-positive solution)
is covered by classical results, and
regularity of the set $\{ u=0\} \cap \{ \nabla u\ne 0\}$ can be obtained via
the implicit function theorem (see Section \ref{non} for higher regularity),
the research focusses on the study of $\partial\{ u>0\}
\cap \partial\{ u<0\} \cap \{ \nabla u=0\}$.\\
In the stationary case --- the two-phase membrane problem ---
the authors proved (\cite{advances} and \cite{imrn}) that 
the free boundary $\partial\{ u>0\} \cup\partial\{ u<0\}$ is in a neighborhood of each branch point,
i.e. a point in the set $\Omega \cap \partial\{ u>0\}
\cap \partial\{ u<0\} \cap \{ \nabla u=0\}$, the
union of (at most) two $C^1$-graphs. Note that the definition of
``branch point'' does not necessarily imply a bifurcation
as that in Figure \ref{branchfig}. 
\newline
We formulate the main result in this paper.
\begin{theorem}\label{main}
Suppose that
$$0<\lambda_{\rm min}\le \inf_{Q_1(0)} \min(\lambda_+,\lambda_-),
\qquad \sup_{Q_1(0)}\max(|\nabla \lambda_+|,|\nabla
\lambda_-|, |\partial_t \lambda_+|,|\partial_t \lambda_-|)<+\infty$$ and that $u$ is a weak solution of
$$\Delta u - \partial_t u = \lambda_+ \chi_{\{u>0\}} \> - \>
\lambda_- \chi_{\{u<0\}}\hbox{ in } Q_1(0)\; ;$$
here $Q_1(0)$ is the parabolic cylinder $(-1,1)\times B_1(0)$.
\newline
Then there are constants $\sigma>0$ and $r_0>0$ such that
\begin{equation}\label{cond} 
u(0)=0\> ,\> |\nabla u(0)|\le \sigma\> ,\> \pardist(0,\{ u>0\})\le \sigma
\> \hbox{ and }\> \pardist(0,\{ u<0\})\le \sigma\end{equation}
imply
$\partial \{ u>0\}\cap Q_{r_0}(0)$ and $\partial
\{ u<0\}\cap Q_{r_0}(0)$ being graphs of
Lipschitz functions (in some space direction) that are continuously differentiable with respect to the space variables.
The constants $\sigma, r_0$, the Lipschitz norms and the modulus of continuity of the spatial normal vectors to these surfaces
depend only on $\inf_{Q_1(0)} \min(\lambda_+,\lambda_-)$, the Lipschitz norms of $\lambda_\pm,$
the supremum norm of $u$ and the space dimension $n$.
\newline
Moreover the regularity above is optimal in the sense that the
graphs are in general not better than Lipschitz.
\end{theorem}
\begin{corollary}
Suppose that
$$0<\lambda_{\rm min}\le \inf_{Q_1(0)} \min(\lambda_+,\lambda_-),
\qquad \sup_{Q_1(0)}\max(|\nabla \lambda_+|,|\nabla
\lambda_-|, |\partial_t \lambda_+|,|\partial_t \lambda_-|)<+\infty$$ and that $u$ is a weak solution of
$$\Delta u -\partial_t u = \lambda_+ \chi_{\{u>0\}} \> - \>
\lambda_- \chi_{\{u<0\}}\hbox{ in } Q_1(0)\; .$$
Then there is a constant $r_0>0$ such that 
if the origin is a branch point, then
$\partial \{ u>0\}\cap Q_{r_0}(0)$ and $\partial
\{ u<0\}\cap Q_{r_0}(0)$ are graphs of
Lipschitz functions (in some space direction)
that are continuously differentiable with respect to the space variables.
The constant $r_0$, the Lipschitz norms and the modulus of continuity of the spatial normal vectors to these surfaces
depend only on $\inf_{Q_1(0)} \min(\lambda_+,\lambda_-)$, the Lipschitz norms of $\lambda_\pm,$
the supremum norm of $u$ and the space dimension $n$.
\end{corollary}
As to the proof we extend the method of \cite{imrn} to the parabolic case.
There is however a difficulty as the time derivative $\partial_t u$
is in general not continuous, so that it is not possible to apply 
directly the comparison principle. We deal with that problem
by a {\em two-stage} proof of directional monotonicity. 
\section{Notation}
Throughout this article $\R^n$ will be equipped with the Euclidean
inner product $x\cdot y$ and the induced norm $\vert x \vert\> ,\>
B_r(x^0)$ will denote the open $n$-dimensional ball of center
$x^0\> ,$ radius $r$ and volume $r^n\> \omega_n\> ,
\> B'_r(0)$ the open $n-1$-dimensional ball of center $0$ and radius
$r\> ,$
and $e_i$ the $i$-th unit vector in $\R^n\> .$
We define $Q_r(t^0,x^0) := (t^0-r^2 , t^0+ r^2)\times B_r(x^0)$ to be 
the cylinder
of radius $r$ and height $2r^2$,
$Q^-_r(t^0,x^0) := (t^0-r^2 , t^0)\times B_r(x^0)$ its ``negative part''
and $Q^+_r(t^0,x^0) := (t^0, t^0+r^2)\times B_r(x^0)$ its ``positive part''.
When omitted, $x^0$ (or $(t^0,x^0)$, respectively) is assumed to be the
origin.
Moreover let $\partial_{\rm par} Q_r(t^0,x^0) := 
(t^0-r^2 , t^0+ r^2)\times \partial B_r(x^0)\cup \{ t^0-r^2 \}\times B_r(x^0)$
denote the parabolic boundary of $Q_r(t^0,x^0)$. 
Let us also introduce the parabolic distance
$\pardist((t,x),A) := \inf_{(s,y)\in A}
\sqrt{\vert x-y\vert^2 + \vert t-s\vert}\> .$
Given a set $A\subset \R^{n+1}\> ,$ we denote its interior by $A^\circ$
and its characteristic function by $\chi_A\> .$
By $\nabla u$ we mean the gradient with respect to the space
variables.
In the text we use the $n$-dimensional Lebesgue-measure
${\mathcal L}^n$ and
the $m$-dimensional Hausdorff measure
${\mathcal H}^m$.
Finally, $\C^{\beta,\mu}:=\H^{\mu,\beta}$
denotes the parabolic
H\"older-space as defined in 
\cite{lady}.
\section{A supremum-mean-value estimate}
In this section we show that at branch points the time derivative
$\partial_t u$, in general a discontinuous function, satisfies a $\sup$-mean-value estimate.
\begin{lemma}\label{mean}
Let $Q^-_{2r}(t^0,x^0)\subset (0,T)\times \Omega$ and let $\lambda_+,\lambda_-$ be non-negative and Lipschitz continuous with respect to the time variable. Then each solution $u$ of (\ref{obst})
satisfies
$$\sup_{Q_1} \vert \partial_t u_{r_k}\vert=\sup_{Q^-_r(t^0,x^0)} |\partial_t u| \le C \left(r^2+\left(r^{-n-2}\int_{Q^-_{2r}(t^0,x^0)}
|\partial_t u|^2\right)^{1\over 2}\right)\; .$$
\end{lemma}
\proof
Using the scaling invariance of the equation
with respect to the scaling
$$u_r(t,x)=r^{-2} u(t^0+r^2t,x^0+rx)$$
we may assume that $r=1/2, t^0=0$ and $x^0=0$.\\
Let $H(t,x,z) = \lambda_+(t,x)\chi_{\{ z>0\}} - \lambda_-(t,x) \chi_{\{ z<0\}}$.
For $$v(t,x) := \partial_t^\tau u(t,x):= {u(t+\tau,x)-u(t)\over \tau}$$
and $\eta\in L^2((-1,1);W^{1,2}(B_1))$ such that $\eta=0$ on 
$(-1,0)\times \partial B_1$,
we calculate
\begin{equation}\label{mean1}
\begin{array}{l}
\int_{-1}^{s} \int_{B_1}(\eta \partial_t v + \nabla v \cdot \nabla \eta)\\
= - \int_{-1}^{s} \int_{B_1}\eta \partial_t^\tau H(t,x,u(t,x))\; , \;  s\in (-1,0)
\; .
\end{array}
\end{equation}
Here
$$\partial_t^\tau H(t,x,u(t,x))
= \lambda_+(t^0+r^2t,x^0+rx) \partial_t^\tau \chi_{\{ u>0\}}
- \lambda_-(t^0+r^2t,x^0+rx) \partial_t^\tau \chi_{\{ u<0\}}$$ $$
+ \chi_{\{ u(t^0+r^2(t+\tau),x^0+rx)>0\}} \partial_t^\tau\lambda_+
- \chi_{\{ u(t^0+r^2(t+\tau),x^0+rx)<0\}} \partial_t^\tau\lambda_-\; .$$
Testing with
$\eta(t,x) := \zeta^2(x)\phi^2(t) \max(v(t,x)-k,0)$ where $k\ge 0$,
$\zeta\in C^{0,1}_0(B_1)$ and $\phi \in C^{0,1}(-1,1)$ such that
$\phi(t)\in [0,1]$ and
$$\phi(t) := \left\{\begin{array}{l}
1, t\ge -1/2\\
0, t\le -1\> ,\end{array}\right.$$
and observing that
$$\max(v(t,x)-k,0) \partial_t^\tau H(t,x,u(t,x))\ge -C_1 r^2 \max(v(t,x)-k,0)$$
we obtain
\begin{equation}\label{supmean}
\sup_{-1 < s < 0} 
\int_{B_1} \phi^2(s) \zeta^2 \max(v(s,\cdot)-k,0)^2
+ \int_{-1}^{s}  \int_{B_1} \phi^2 \zeta^2 |\nabla \max(v-k,0)|^2
\end{equation}
$$\le C_2 \int_{-1}^{s}  \int_{B_1} [
\max(v-k,0)^2 (\phi^2 |\nabla \zeta|^2 | + 
\phi |\partial_t \phi| \zeta^2)+ r^2 \phi^2\zeta^2\max(v-k,0)]\; .$$
From the proof of \cite[Theorem 4.7]{lieberman} we infer that
\begin{equation}\label{time_1}\sup_{Q_{1/2}^-} v
\le C_3  \left(r^2+\left(\int_{Q_1^-} v^2\right)^{1\over 2}\right)\; .\end{equation}
Testing with
$\eta(t,x) := \zeta^2(x)\phi^2(t) \max(-v(t,x)-k,0)$ where $k\ge 0$,
we obtain in a similar way that
\begin{equation}\label{time_2}
\sup_{Q_{1/2}^-} (-v)
\le C_3  \left(r^2+\int_{Q_1^-} v^2\right)^{1\over 2}\; .
\end{equation}
Letting $\tau\to 0$ and scaling back we obtain the statement.
\qed
\section{Non-degeneracy and regularity of the solution}
\begin{lemma}[Non-Degeneracy]
\label{ndeg}
For every $Q_{2r}(t^0,x^0)\subset (0,T)\times \Omega$ the following holds:
\[ \textrm{1) If } (t^0,x^0)\in \partial \{ u>0\}, \textrm{ then }\sup_{Q_r^-(t^0,x^0)} u\; \ge \;
{1\over {8n}} \inf_{Q_{r}(t^0,x^0)}\lambda_+ \> r^2\; .\]
\[ \textrm{2) If } (t^0,x^0)\in \partial \{ u<0\}, \textrm{ then } \inf_{Q_r^-(t^0,x^0)} u\; \le \;
-{1\over {8n}}\inf_{Q_{r}(t^0,x^0)}\lambda_- \> r^2\; .
\]
\end{lemma}
\proof 
We choose a sequence
$\{u>0\}\ni (t^m,x^m) \to (t^0,x^0)$ as $m\to \infty\> .$
Supposing that $\sup_{Q_r^-(t^m,x^m)} u \> \le \>
{1\over {8n}} \inf_{Q_{r}(t^0,x^0)}\lambda_+ \> r^2\> ,$ the comparison principle
yields that $u(t,x)\le v(t,x) := ({t^m-t\over 2}+{1\over {8n}}{\vert x-x_m\vert}^2)\>\inf_{Q_{r}(t^0,x^0)}\lambda_+
$ in $Q_r^-(t^m,x^m)\> ,$ a contradiction to the fact that $u(t^m,x^m)>0\> .$\\
The estimate for $\inf_{Q_r^-(t^0,x^0)} u$ is obtained the same
way, replacing $u$ by $-u$ and $\lambda_+$ by $\lambda_-\> .$ \qed
\begin{lemma}\label{bounded}
Let $\lambda_+,\lambda_-\in C^{0,1}_{\rm loc}((0,T)\times \Omega)$. Then each solution $u$ of (\ref{obst})
satisfies the following:\\
1) $\partial_t u\in L^\infty_{\rm loc}((0,T)\times \Omega)$.\\
2) $\partial_t\nabla u\in L^2_{\rm loc}((0,T)\times \Omega)$.
\end{lemma}
\proof
1) follows from Lemma \ref{mean}.\\
2) follows from (\ref{supmean}) with $k=0$ and from the analogous estimate for $\max(-v,0)$.
\qed
\begin{corollary}\label{ndeg2}
\label{ndeg2}
For every $Q_{2r}(t^0,x^0)\subset (0,T)\times \Omega$,
there exists a constant $c_0>0$ depending only on $n$ and 
$\Vert\partial_t u\Vert_{L^\infty(Q_r(t^0,x^0))}$
such that
\[ u \ge 0 \textrm{ in } Q_r^-(t^0,x^0) \textrm{ implies } u\ge 0 \textrm{ in } Q_{c_0r}(t^0,x^0)\; ,\textrm{ and}\]
\[ u \le 0 \textrm{ in } Q_r^-(t^0,x^0) \textrm{ implies } u\le 0 \textrm{ in } Q_{c_0r}(t^0,x^0)\; .\]
\end{corollary}
\proof
Suppose towards a contradiction that
$u(t^1,x^1)<0$ for some $(t^1,x^1)\in Q_{c_0r}^+(t^0,x^0)$. Then there is a point
$(t^2,x^2)\in \partial\{ u<0\}\cap \overline{Q_{c_0r}^+(t^0,x^0)}$. Applying Lemma \ref{ndeg}
at $(t^2,x^2)$ with respect to the cylinder $Q_{(1-c_0)r}(t^2,x^2)$ yields a contradiction to Lemma \ref{bounded} 1) provided
that $c_0$ has been chosen small enough.\\
The second estimate is proved in the same fashion.\qed
\begin{proposition}\label{regular}
 Let $\lambda_+,\lambda_-\in C^{0,1}_{\rm loc}((0,T)\times \Omega)$. Then each solution $u$ of (\ref{obst})
satisfies $\nabla u\in {\bf C}^{1/2,1}_{\rm loc}((0,T)\times \Omega)$, that is, the gradient is Lipschitz continuous with respect to the space variables and H\"older continuous with exponent $1/2$ with respect to the time variable.
\end{proposition}
\proof
Let us first show that for any $e\in \partial B_1$,
$(\Delta-\partial_t)(\max(\partial_e u,0))\ge -C$ and
$(\Delta-\partial_t)(\max(-\partial_e u,0))\ge -C$ in $\Omega$.
We give a formal proof that can be made rigorous translating
everything into a weak formulation. In $\{ \partial_e u > 0\}$,
$$(\Delta-\partial_t)(\partial_e u)
$$ $$= {\partial_e u \over {|\nabla u|}} (\lambda_+ {\mathcal H}^{n-1}\lfloor (\{ \nabla u\ne 0\} \cap \partial\{ u>0\})
+ \lambda_- {\mathcal H}^{n-1}\lfloor (\{ \nabla u\ne 0\} \cap \partial\{ u<0\}))$$
$$+ \partial_e \lambda_+ \chi_{\{ u>0\}} - \partial_e \lambda_- \chi_{\{ u<0\}}\; 
\ge \; -C\; .$$
As $\partial_e u$ is continuous,
we obtain $(\Delta-\partial_t)(\max(\partial_e u,0))\ge -C$.\\
Considering $-e$ instead of $e$ we obtain also
$(\Delta-\partial_t)(\max(-\partial_e u,0))\ge -C$.
But then the ``almost monotonicity formula'' Theorem I of \cite{edquist} applies
and we proceed as follows:
at each point $(t^0,x^0)\in \{ u\ne 0\}\cap \{ \nabla u=0\}$,
we obtain from the almost monotonicity formula that
$\nabla \partial_e u$ is
bounded at $(t^0,x^0)$ by a locally uniform constant.
\\
At each point $(t^0,x^0)\in \{ u\ne 0\}\cap \{ \nabla u\ne 0\}$,
we obtain in a similar way that
{\em for every $e\bot \nabla u(t^0,x^0)$},
$|\nabla \partial_e u(t^0,x^0)|$ is
bounded by a locally uniform constant.
Let $e_1 = {\nabla u(t^0,x^0)\over {|\nabla u(t^0,x^0)|}}$.
Then $-\partial_{11} u(t^0,x^0) = -\lambda_+ \chi_{\{ u(t^0,x^0)>0\}}
+ \lambda_- \chi_{\{ u(t^0,x^0)<0\}} - \partial_t u(t^0,x^0) + \sum_{j=2}^n \partial_{jj} u(t^0,x^0)$
is by Lemma \ref{bounded} bounded by a locally uniform constant.
\qed
\begin{corollary}\label{density}
${\mathcal L}^{n+1}(\partial \{ u > 0\}\cup \partial \{ u< 0\})=0$
\end{corollary}
\proof First, we obtain from Lemma \ref{ndeg}, Lemma \ref{bounded}
and Proposition \ref{regular} that there exists a locally uniform
constant $c>0$ such that for $Q_{2r}(s,y)\subset (0,T)\times \Omega$,
$$\frac{{\mathcal L}^{n+1}(Q_r(s,y)\cap \{ u>0\})}{{\mathcal L}^{n+1}(Q_r)} \ge c >0\textrm{ if }(s,y)\in \partial \{ u> 0\}$$
$$\textrm{and }\frac{{\mathcal L}^{n+1}(Q_r(s,y)\cap \{ u<0\})}{{\mathcal L}^{n+1}(Q_r)} \ge c >0\textrm{ if }(s,y)\in \partial \{ u< 0\}\; .$$
Since $\chi_{\{ u>0\}} * \chi_{Q_r}/{\mathcal L}^{n+1}(Q_r)\to \chi_{\{ u>0\}}$ in $L^1_{\rm loc}((0,T)\times \Omega)$ as $r\to 0$
and the analogous fact holds for $\chi_{\{ u<0\}}$, we obtain
that $\chi_{\{ u>0\}}\ge  c >0$ ${\mathcal L}^{n+1}$-a.e. on $\partial \{ u> 0\}$
and $\chi_{\{ u<0\}}\ge  c>0$ ${\mathcal L}^{n+1}$-a.e. on $\partial \{ u< 0\}$. Thus 
${\mathcal L}^{n+1}(\partial \{ u > 0\}\cup \partial \{ u < 0\})=0$.
\qed
\section{Vanishing time derivative}
As a corollary of Lemma \ref{mean} we obtain now that
at points at which the blow-up limit depends only
on the space variables, the time derivative
$\partial_t u$ -- in general a discontinuous function -- 
attains the limit $0$.
\begin{corollary}\label{zerolim}
Let $Q_{2r}(t^0,x^0)\subset (0,T)\times \Omega$
and suppose that for a sequence of solutions $u_k$ in $(0,T)\times \Omega$
$$u_{r_k}(t,x)={r_k}^{-2} u_k(t^k+r_k^2t,x^k+r_kx)\to u_0(x)
\textrm{ in } L^1_{\rm loc}(\R^{n+1})\textrm{ as }r_k\to 0\; .$$
Then
$$\sup_{Q_{r_k}(t^k,x^k)} |\partial_t u_k| \to 0$$
as $r_k\to 0$.
\end{corollary}
\proof
The statement follows from Lemma \ref{mean}
and the fact that $\partial_t u_{r_k}$ converges to $0$ in $L^2_{\rm loc}(\R^{n+1})$
as $r_k\to 0$.
The $L^2$-convergence in turn may be shown as follows:
as $\partial_t u_k$ is by Lemma \ref{mean} bounded in $L^\infty(Q_{r}(t^0,x^0))$, it is sufficient to prove
a.e. convergence.
For $(s,y)\in \{ u_0 =0\}^0$ we obtain from Lemma \ref{ndeg} that
$u_{r_k}=0$ in $Q_\delta(s,y)$ for some $\delta>0$ and large $k$.
For $(s,y)\in \{ u_0 > 0\}\cup \{ u_0 < 0\}$, $u_{r_k}$ converges in $C^1(Q_\delta(s,y))$
for some $\delta>0$ as $k\to\infty$.
Moreover we know from Corollary \ref{density} that ${\mathcal L}^{n+1}(\partial \{ u_0 > 0\}\cup \partial \{ u_0< 0\})=0$.
It follows that $\partial_t u_{r_k}$ converges ${\mathcal L}^{n+1}$-a.e. to $\partial_t u_0$.
\qed
\section{Directional monotonicity}\label{dirmon}
In a first stage, we show that 
if the solution is close to the one-dimensional solution
\begin{equation}\label{one-dim-sol}
h(x) :={\lambda_+(0)\over 2} \max(x_1,0)^2\> -{\lambda_-(0)\over 2}
\min(x_1,0)^2 \ .
\end{equation}
then it is increasing in a cone of {\em spatial} directions.
Later on we will extend the result to a cone of {\em tempo-spatial} directions.
\begin{proposition}\label{directional-monoton}
Let $0<\lambda_{\rm min}\le \inf_{Q_1(0)}
\min(\lambda_+,\lambda_-)$, $h$ as in (\ref{one-dim-sol}), and let
$\varepsilon\in (0,1)$. Then each solution $u$ of (\ref{obst}) in
$Q_1(0)$ such that
$$\dist_{L^\infty((-1,1);W^{1,\infty}(B_1))}(u,h)
\le \delta := {\lambda_{\rm min}\varepsilon\over {48 n}}$$ and
$$\sup_{Q_1(0)}\max(|\nabla \lambda_+|,|\nabla \lambda_-|)\le
\delta$$ satisfies $\varepsilon^{-1}
\partial_e u - |u| \ge 0$ in $Q_{1/2}(0)$ for every $e\in \partial
B_1(0)$ such that $e_1\ge \varepsilon$;
here $e_1$ denotes the first component of the vector $e$.
\end{proposition}
\proof
First note that $\varepsilon^{-1} \partial_e h - |h|\ge 0$ in $Q_2(0)$.
It follows that
\begin{equation}\label{approx}
\varepsilon^{-1} \partial_e u - |u|\ge -3\delta \varepsilon^{-1} \textrm{ in } Q_1(0)
\end{equation}
provided that $\dist_{L^\infty((-1,1);W^{1,\infty}(B_1))}(u,h) \le\delta$. Suppose now
towards a contradiction that the statement is not true. Then there
exist $\lambda_+,\lambda_-\in (\lambda_{\rm min},+\infty), (t^*,x^*) \in
Q_{1/2}(0), e^*, $ and a solution $u$ of (\ref{obst}) in $Q_1(0)$
such that $\dist_{L^\infty((-1,1);W^{1,\infty}(B_1))}(u,h) \le\delta$,
$$\sup_{Q_1(0)}\max(|\nabla \lambda_+|,|\nabla \lambda_-|)\le
\delta,$$ $e^*_1 \ge \varepsilon$ and
$\varepsilon^{-1}\partial_{e^*} u(t^*,x^*) - |u(t^*,x^*)|<0.$ For the
positive constant $c$ to be defined later the functions $v :=
\varepsilon^{-1}\partial_{e^*} u - |u|$ and $w:=
\varepsilon^{-1}\partial_{e^*} u - |u|+c |x-x^*|^2 - c (t-t^*)$ satisfy then
the
following: in the set $D := Q_1(0)\cap \{ v<0\}\cap\{ t<t^*\}$,
\[ \Delta w \> -\> \partial_t w\le 2nc + c -{\lambda_+} \chi_{\{ u>0\}}-{\lambda_-} \chi_{\{ u<0\}}\]\[
+\varepsilon^{-1} (\lambda_+ +\lambda_-) \nu_x\cdot e^* {\mathcal
H}^{n-1}\lfloor (\{ u=0\}\cap \{ \nabla u\ne 0\}) \]\[+\varepsilon^{-1}
( \chi_{\{ u>0\}}\partial_{e^*} \lambda_+ - \chi_{\{
u<0\}}\partial_{e^*} \lambda_-)
\]
where $\nu_x = {\nabla u \over {|\nabla u|}}$. As
\[\nu_x\cdot e^*< 0 \> \hbox{ on }\> \{ u=0\}\cap \{ v<0\} = \{ u=0\}\cap \{
\varepsilon^{-1}\partial_{e^*} u <0\}\; ,\] we obtain by the definition of $\delta$ that
$w$ is supercaloric in $D$ provided that $c$ has been chosen
accordingly, say $c:= \lambda_{\rm min}/(4n)$. It follows that the
negative infimum of $w$ is attained on
\[ \partial_{\rm par} D \subset (\partial_{\rm par} Q_1(0)\cap \{ t\le t^*\}) \cup
\left( Q_1(0)\cap \partial \{ v<0\}\right)
\; .\]
Consequently it is attained on $\{ t\le t^*\}\cap \partial_{\rm par} Q_1(0)$,
say at the point $(\bar t,\bar x) \in \{ t \le t^*\}\cap \partial_{\rm par} Q_1(0)$.
Since $\pardist((\bar t,\bar x),(t^*,x^*))\ge 1/2$, we obtain that
\[ \varepsilon^{-1}\partial_{e^*} u(\bar t,\bar x) - |u(\bar t,\bar x)|=
v(\bar t,\bar x) = w(\bar t,\bar x) - c |x^* - \bar x|^2 + c (\bar t-t^*) 
\]\[ < -c/4=-\lambda_{\rm min}/(16n)\; .\] But this contradicts
(\ref{approx}) in view of $\delta = {\lambda_{\rm
min}\varepsilon\over {48 n}}$.\qed
\section{The set of non-vanishing gradient}\label{nonvanish}
In the sequel we are going to need higher regularity of
the level set $\{ u=0\}\cap \{\nabla u\ne 0\}$.
Higher regularity can be obtained in a standard way using the
von Mises transform:
\begin{lemma}\label{cinf}
The set $\{ u=0\}\cap \{\nabla u\ne 0\}$ is locally in $(0,T)\times \Omega$ a $C^1$-surface
and $\partial_t u$ is continuous on that surface.
\end{lemma}
\proof
Let $(t^0,x^0)\in \{ u=0\}\cap \{\nabla u\ne 0\}$.
We may assume that $\nabla u(t^0,x^0)=\partial_1 u(t^0,x^0)$ and
that in $Q_\delta(t^0,x^0)$, 
$u$ is strictly increasing in the $x_1$-direction and
$\{ u=0\}$ is the graph
of a function, say $x_1=g(t,x')$ for $(t,x)\in Q_\delta(t^0,x^0)$,
where $g\in C^0((t^0-\delta^2,t^0+\delta^2);C^1(B'_\delta(x^0)))$.
It is sufficient to prove that
$g\in C^1(Q_{\delta/2}(t^0,x^0))$.
To do so, we use von Mises variables, i.e.
$$y=u(t,x_1,x') \textrm{ and } x_1=v(t,y,x')\; .$$
A calculation assures that
$$\left(\frac{-1-|\nabla' v|^2}{(\partial_y v)^3} \partial_{yy} v\right)
-\frac{\Delta' v}{\partial_y v} + 2 \frac{\nabla' v \cdot \nabla' \partial_y v}{(\partial_y v)^2}
+ \frac{\partial_t v}{\partial_y v} $$ $$= \left\{\begin{array}{l}
\lambda_+(t,v(t,y,x'),x') , \> y>0\\
-\lambda_-(t,v(t,y,x'),x') , \> y<0\end{array}\right. \; .$$
Thus
$$\partial_t v - a_{ij}(\nabla v) \partial_{ij} v =
f(t,y,x')\partial_y v:=
\left\{\begin{array}{l}
-\lambda_+(t,v(t,y,x'),x') \partial_y v, \> y>0\\
 \lambda_-(t,v(t,y,x'),x') \partial_y v, \> y<0\> .\end{array}\right.$$
Provided that $\delta$ has been chosen small enough,
$|\nabla' v|\le 1/2$, $0< \partial_y v \le 1/2$ and
the above equation is uniformly parabolic.
Moreover,
$$\partial_t \partial^h_t v - a_{ij}(\nabla v) \partial_{ij} \partial^h_t v 
- \frac{\partial a_{ij}(z_h)}
{ \partial p_k}\partial_{ij} v(t+h,y,x') \partial_k \partial^h_t v$$
$$=
f(t,y,x')\partial_y \partial^h_t v + \partial_y v(t+h,y,x')\partial^h_t f(t,y,x')$$
where $z_h=\theta(t,y,x') \nabla v(t+h,y,x') + (1-\theta(t,y,x'))  \nabla v(t,y,x')$
and $\theta(t,y,x')\in [0,1]$.
Calculating 
$$\partial_1 v = 1/\partial_1 u, \partial_t v = -\partial_t u/\partial_1 u, \partial_i v = -\partial_i u/\partial_1 u
\textrm{ for } 2\le i \le n,$$
$$\partial_{ij} v + \partial_{i1} v \partial_j u + \partial_{1j}v\partial_i u +\partial_1 v \partial_{ij} u+\partial_{11}v\partial_i u\partial_j u=0
\textrm{ for } 2\le i,j \le n
$$
$$\partial_{i1} v \partial_1 u + \partial_{11} v \partial_1 u \partial_i u + \partial_1 v \partial_{i1} u =0,
\partial_{11} v = -\partial_{11} u/(\partial_1 u)^3$$
shows that all spatial second derivatives of $v$ 
and $\partial_t v$ are bounded.
Thus $f(t,y,x'), \partial^h_t f(t,y,x')$ and $\frac{\partial a_{ij}(z_h)}
{ \partial p_k}\partial_{ij} v(t+h,y,x')$
are bounded uniformly in $h$, and we obtain from \cite{krylov}
that $\partial^h_t v$ is uniformly H\"older continuous with respect to $h$
and that $\partial_t u$ is H\"older continuous in $Q_{\delta/2}(t^0,x^0)$.
\qed
\section{Global solutions}
In this section we extend our characterization of elliptic global solutions
\cite[Theorem 4.3]{global} to the parabolic case.
We are going to need the following version of the Caffarelli-Kenig
monotonicity formula of \cite{caffarellikenig}:
\begin{theorem}\label{mono}
Let
$$\Phi(r,w) := {1\over {r^4}} I(r,\max(w,0))I(r,\max(-w,0))$$
where
$$I(r,v) := \int_{-r^2}^0 \int_{\R^n} |\nabla v|^2 G(t,x)$$
and $G$ is the backwards heat kernel
$$G(t,x) = (4\pi (-t))^{n/2} \exp({|x|^2\over {4t}})\; .$$
If $\max(w,0)$ and $\max(-w,0)$ are continuous subcaloric functions, then
$r\mapsto \Phi(r,w)$ is non-decreasing, and
$\Phi(\sigma,w)=\Phi(\rho,w)$ for some $0<\rho<\sigma$ implies that
either
\\
(A) $\nabla \max(w,0)=0$ in $-\sigma^2<t<0$ or $\nabla \max(-w,0)=0$ in $-\sigma^2<t<0$.\\
or
\\
(B) $\max(w,0)(\partial_t - \Delta)\max(w,0)=0$
and $\max(-w,0)(\partial_t - \Delta)\max(-w,0)=0$
in $-\sigma^2<t<0$
in the sense of measures.
\end{theorem}
\proof
For $v := \max(w,0)$ (or $v := \max(-w,0)$, respectively) we calculate
$$I(r,v)=
- {1\over 2} \int_{-r^2}^0 \int_{\R^n}  G(t,x)(\partial_t - \Delta)v^2
+ \int_{-r^2}^0 \int_{\R^n}  G(t,x)v(\partial_t - \Delta)v\; ,$$
$$I'(r,v)\ge 2r\int_{\R^n} |\nabla v|^2 G(-r^2,x)\; .$$
In what follows we assume that $I(r,v)\ne 0$.
It follows that
$$\frac{I'(r,v)}{I(r,v)} \ge 4r\frac{\int_{\R^n} |\nabla v(-r^2,x)|^2 G(-r^2,x)}{\int_{\R^n}  v^2(-r^2,x)G(-r^2,x)}\; .$$
In the case $I'(r,v)\ne 0$ the inequality is strict unless $\int_{-r^2}^0 \int_{\R^n}v(\partial_t - \Delta)v=0$.
Consequently $\Phi(r,w)=0$, or else
$$\frac{\Phi'(r,w)}{\Phi(r,w)} \ge {4\over r}
\Bigg[-1\;+ \; r^2 \frac{\int_{\R^n} |\nabla \max(w,0)|^2 G(-r^2,x)}{\int_{\R^n} \max(w,0)^2G(-r^2,x)}$$
$$+ r^2 \frac{\int_{\R^n} |\nabla \max(-w,0)|^2 G(-r^2,x)}{\int_{\R^n} \max(-w,0)^2G(-r^2,x)}\Bigg]\; ,$$
where the inequality is strict unless both $\int_{-r^2}^0 \int_{\R^n}\max(w,0)(\partial_t - \Delta)\max(w,0)=0$
and $\int_{-r^2}^0 \int_{\R^n}\max(-w,0)(\partial_t - \Delta)\max(-w,0)=0$.
Moreover, by \cite[Corollary 2.4.6]{caffarellikenig}, the right-hand side is non-negative.
\qed
\begin{lemma}\label{forward}
Let $v^1,v^2$ be solutions of (\ref{obst}) in $\R^{n+1}$ with
such that
$v^1=v^2$ in $\{ t<0\}$ and $v^1,v^2$ have polynomial growth with respect to the space variables.
Then $v^1=v^2$ in $\R^{n+1}$.
\end{lemma}
\proof
Multiplying the difference of the two equations by $(v^1-v^2)W$ where $W(t,x)=G(t-T,x)$
and
integrating, we obtain for each $0<T<+\infty$, $0<S<T$ and $H$ defined in Lemma \ref{mean} that $0\; =$
$$ \int_0^S \int_{\R^n}  W [|\nabla (v^1-v^2)|^2 \> + \> (H(v^1)-H(v^2))(v^1-v^2)]
\; - \; {1\over 2}\int_0^S \int_{\R^n} (v^1-v^2)^2 \partial_t W$$ $$
+ \; {1\over 2} \int_{\R^n} W(S)  (v^1(S)-v^2(S))^2
\; + \; \int_0^S \int_{\R^n} (v^1-v^2) \nabla W \cdot \nabla (v^1-v^2)
$$
$$\ge {1\over 2} \int_{\R^n} W(S)  (v^1(S)-v^2(S))^2\; + \; {1\over 2}\int_0^S \int_{\R^n} (v^1-v^2)^2 [-\partial_t W-\Delta W]$$
$$= \; {1\over 2} \int_{\R^n} W(S)  (v^1(S)-v^2(S))^2\; .$$\qed
\begin{lemma}\label{backself}
Assume that $w$ is a backward self-similar solution with constant coefficients $\lambda_+,\lambda_-$,
i.e.
$$ w(\theta^2 t,\theta x)= \theta^2 w(t,x) \textrm{ for all }
\theta\ge 0, t<0 \textrm{ and } x\in \R^n\; .$$
Then $\nabla w=0$ on $\{ w=0\}$.
\end{lemma}
\proof
First, the self-similarity implies that
\begin{equation}\label{homog}
\partial_e w(\lambda^2t,\lambda x)=\lambda \partial_e w(t,x)
\textrm{ for all }e\in \partial B_1, \lambda \ge 0, t<0 \textrm{ and }x\in \R^n.
\end{equation}
Consequently the function $r\mapsto \Phi(r, \partial_e w)$ of 
the monotonicity formula Theorem \ref{mono}
is constant in $(0,+\infty)$, implying by Theorem \ref{mono}
that 
either
\\
(A) $\nabla \max(\partial_e w,0)=0$ in $\{ t<0\}$ or $\nabla \max(-\partial_e w,0)=0$ in $\{ t<0\}$.
\\
or\\
(B) $\max(\partial_e w,0)(\partial_t - \Delta)\max(\partial_e w,0)=0$ in $\{ t<0\}$ 
and $\max(-\partial_e w,0)(\partial_t - \Delta)\max(-\partial_e w,0)=0$ in $\{ t<0\}$ 
in the sense of measures.\\
Suppose now towards a contradiction that 
there is a point $(t^1,x^1)\in \{ t<0\} \cap \{ w=0\} \cap \{ \nabla w \ne 0\}$
and denote $\nu={\nabla w\over {|\nabla w|}}$,
$\nu^0= {\nabla w(t^1,x^1)\over {|\nabla w(t^1,x^1)|}}$ and let $Q_\kappa(t^1,x^1)$ such
that $\partial_{\nu^0} w >0$ in $Q_\kappa(t^1,x^1)$ and
$\{ w=0\}\cap Q_\kappa(t^1,x^1)$ is a $C^1$-surface.
In the case $\nu^0\cdot e\ne 0$,
$$ |(\partial_t - \Delta)\partial_e w|(Q_\kappa(t^1,x^1))
= 
|\lambda_++\lambda_-|
\int_{t^1-\kappa^2}^{t^1+\kappa^2}
\int_{B_\kappa(x^1)\cap \{ w(t)=0\} } |e\cdot \nu| \> d {\mathcal H}^{n-1}\> dt
\ne 0\; .$$
Thus (A) holds.
From (\ref{homog}) we infer that $\partial_e w\ge 0$ in $\{ t<0\}$ if $e\cdot \nu^0>0$ and
$\partial_e w\le 0$ in $\{ t<0\}$ if $e\cdot \nu^0<0$.
Hence $\partial_e w= 0$ in $\{ t<0\}$ for all $e\bot \nu^0$. 
As in \cite[p. 844]{cps} we may
write
$$ w(t,x)= -t f({x_n\over{\sqrt{-t}}})$$
and calculate the $2$-parameter family of solutions
of the ODE which $f(\xi)=w(-1,\xi)$ satisfies in $(0,+\infty)$,
$$f(\xi)=\lambda_++C_1 (\xi^2-2)$$ $$+\; C_2\left( -2\xi e^{\xi^2/4} + (\xi^2-2) \int_0^{\xi} e^{s^2/4}\> ds\right)\; \textrm{ in } \{ f>0\}$$
and
$$f(\xi)=-\lambda_-+C_3 (\xi^2-2)$$ $$+\; C_4\left( -2\xi e^{\xi^2/4} + (\xi^2-2) \int_0^{\xi} e^{s^2/4}\> ds\right)\; \textrm{ in } \{ f<0\}\; .$$
As $w$ has polynomial growth towards infinity
we conclude that
$0=C_2=C_4$ and that
$$f(\xi)=\lambda_++C_1 (\xi^2-2) \textrm{ in } \{ f>0\}$$
and
$$f(\xi)=-\lambda_-+C_3(\xi^2-2) \textrm{ in } \{ f<0\}\, .$$
If $f(a)=0$ and $f'(a)\ne 0$ for some $a\in \R$ then
$C_1=C_3=-\lambda_+/(a^2-2)=\lambda_-/(a^2-2)$, a contradiction.
Therefore $f(a)=0$ implies $f'(a)=0$ and $a=0$.
It follows that
$\nabla w=0$ on $\{ w=0\}$.\qed
\begin{theorem}\label{global-solution}
Let $w$ be a global solution 
with constant coefficients $\lambda_+,\lambda_-$
such that $\partial_t w$ and $D^2 w$ are bounded,
and suppose that
the origin (in time-space) is a branch point of $w$.
Then after rotation 
$$w(t,x)=w^*(t,x) := \lambda_+ {\max(x_n,0)^2/2} - \lambda_- {\max(-x_n,0)^2/2}
\textrm{ for } (t,x)\in \R^{n+1} \; .$$
\end{theorem}
\noindent {\sl Proof.}\\
{\bf Step 1:}  Let us first assume that $w$ is a backward self-similar solution.
By Lemma \ref{backself}  $\nabla w=0$ on $w=0$.
But then $z_1:= \max(w,0)$ and $z_2:= \max(-w,0)$ are in $\{ t \le 0\}$ non-negative backward self-similar solutions. Concerning those, it has been
shown in \cite[Lemma 6.3]{cps} and \cite[Theorem 8.1]{cps} that
either $z_j$ is a half-plane solution of the form
$z_j(t,x) = \lambda_\pm/2 \max(x\cdot e,0)^2$ for some
$e\in \partial B_1$, or
$z_j(t,x)=- a_0 t + \sum_{i=1}^n a_i x_i^2$ with non-negative 
constants $a_i, 0\le i\le n$.
In the latter case the symmetry of $z_j$ implies that
$z_k=0$ in $\{ t<0\}$ for $k\ne j$, and by Corollary
\ref{ndeg2} the origin cannot be a branch point. 
\\
It follows that after rotation
$$w(t,x)= w^*(t,x)\textrm{ for } t<0 \; .$$ 
\\
{\bf Step 2:}
In the case of a general solution $w$ as in the statement of our theorem,
we consider the blow-up up $w_0$ of $w$ at the origin and the blow-down
$w_\infty$. By the non-degeneracy Lemma \ref{ndeg} and 
\cite[Theorem 4.1]{siam},
both $w_0$ and $w_\infty$ satisfy 
the assumptions of Step 1.
Thus both $w_0$ and $w_\infty$ are after rotation
of the form
$\lambda_+ {\max(x_n,0)^2/2} - \lambda_- {\max(-x_n,0)^2/2}$ for $t<0$,
and the monotonicity formula
\cite{siam} implies 
that $w$ is backward self-similar. But then it follows from Step 1
that  
after rotation 
$$w(t,x)=w^*(t,x)
\textrm{ for } t<0 \; .$$
Last, we apply Lemma \ref{forward} to obtain the same
for $t\ge 0$.
\qed
\section{Uniform closeness to $h$}
We are now ready to prove uniform closeness of the scaled solution
to the global solution $h$ of (\ref{one-dim-sol}), assuming that we are in the setting
of Theorem \ref{main}.
\begin{lemma}\label{uniform-closeness}
Let $u$ be a solution of (\ref{obst}) in $Q_1(0)$. Then, given
$\delta >0$, there are constants $r_\delta
>0,\sigma_\delta>0$ (depending only on $\inf_{Q_1(0)} \min(\lambda_+,\lambda_-)$,
the Lipschitz norms of $\lambda_\pm,$ the supremum norm of $u$ and
the space dimension $n$) such that the following holds:
\newline
If $r\in (0,r_\delta]\> , \> u(s,y)=0\> ,\> |\nabla u(s,y)|\leq \sigma_\delta r$,
$\pardist((s,y),\{ u>0\})\le \sigma_\delta r$ and $\pardist((s,y),\{ u<0\})\le \sigma_\delta r$
for some $(s,y) \in
Q_{1/2}(0)$ then in
$Q_{r}(s,y)$,
 the solution
$u(s+\cdot,y+\cdot )$ is $\delta r^2$-close to a rotated version $\tilde h$ of the
one-dimensional solution $h$ defined in (\ref{one-dim-sol}), more precisely
$$
r^{-2}\sup\limits_{Q_{r}(0)}|u(s+\cdot,y+\cdot)-\tilde h|+
r^{-1}\sup\limits_{Q_{r}(0)}|\nabla u(s+\cdot,y+\cdot)-\nabla {\tilde h}|+
\sup\limits_{Q_{r}(0)}|\partial_t u(s+\cdot,y+\cdot)|$$
$$\leq \delta.
$$
\end{lemma}
\proof Suppose towards a contradiction that the statement of the lemma fails. Then for some $\delta
>0$ there exist $\sigma_j\to 0, r_j \to 0$, $(s^j,y^j) \to (s^0,y^0)\in \overline{Q_{1/2}}$, a sequence $u_j$
 of solutions such that
$(s^j,y^j)\in Q_{1/2}(0)$,
$u_j(s^j,y^j)=0$,
$|\nabla u_j(s^j,y^j)|\leq \sigma_j r_j$,
$\>\pardist((s^j,y^j),\{ u_j>0\})\le \sigma_jr_j$, $\>\pardist((s^j,y^j),\{ u_j<0\})\le \sigma_jr_j$
 and
$$
r^{-2}_{j}\sup\limits_{Q_1(0)}| u_j(s^j+r_j^2 \cdot,y^j+r_j\cdot)-\tilde h
(r_j\cdot)| \> + \> r^{-1}_{j }\sup\limits_{Q_1(0)}| \nabla
u_j(s^j+r_j^2 \cdot,y^j+r_j\cdot)-\nabla \tilde h (r_j\cdot)|
$$
$$ + \> \sup\limits_{Q_1(0)}| \partial_t u_j(s^j+r_j^2 \cdot,y^j+r_j\cdot)|
 > \delta
$$
for all possible rotations $\tilde h$ of $h$.
\\
We may define
$$
U_j(x):=\frac{u_j(r_j^2 t + s^j, r_j x + y^j)}{r_j^{2}}
$$
and arrive at
\begin{equation}\label{hdiff}
\Vert U_j-\tilde h\Vert_{W^{1,\infty}(Q_{1})} > \delta ,
\end{equation}
for all possible rotations $\tilde h$ of $h$.
\newline
Observe that $U_j$ is a solution of (\ref{obst}) in $Q_1$ with
respect to the scaled coefficients $\lambda_+(r_j^2 t+s^j,r_j x + y^j)$ and
$\lambda_-(r_j^2 t+s^j, r_j x + y^j)$. Since $U_j(0)=0$, $\>|\nabla U_j(0)|\leq
\sigma_j$,$\>\pardist(0,\{ U_j>0\})\le \sigma_j$,$\>\pardist(0,\{ U_j<0\})\le \sigma_j$ and the derivatives $D^2 U_j, \partial_t U_j$ are uniformly
bounded, we obtain by standard compactness arguments a global limit
solution $U_0$ of (\ref{obst}) in $\R^n$ with respect to
$\lambda_+(s^0,y^0)$ and $\lambda_-(s^0,y^0)$ which
satisfies $0\in \partial\{ U_0>0\}\cap \partial\{ U_0<0\}\cap \{ \nabla U_0=0\}$.
By Theorem \ref{global-solution},
$U_0=\tilde h$ where $\tilde h$ is a rotated version of $h$.
Thus $U_j$ and $\nabla U_j$ converge in $Q_1$ uniformly to $\tilde h$
and $\nabla \tilde h$, respectively, and by Corollary \ref{zerolim}
$\partial_t U_j \to 0$ in $L^\infty(Q_1)$ as  $j\to\infty$.
We obtain a contradiction to (\ref{hdiff}).\qed
\section{Continuity of the time derivative}
Assuming once more that we are in the setting of Theorem \ref{main},
we show in the present section that the time derivative 
of the solution is {\em continuous} in a suitable neighborhood of the origin.
\begin{proposition}\label{time-indep}
Let $u$ be a solution of (\ref{obst}) in $Q_1$.
Then
there are positive constants $\tilde r$ and $\tilde \sigma$
(depending on $\inf_{Q_{1}} \min(\lambda_+,\lambda_-)$,
the Lipschitz norms of $\lambda_\pm,$ the supremum norm of $u$ and
the space dimension $n$)
such that
the following holds.
If $u(0)=0\> ,\> |\nabla u(0)|\leq \tilde \sigma\tilde r$,
$\pardist(0,\{ u>0\})\le \tilde \sigma\tilde r$ and $\pardist(0,\{ u<0\})\le \tilde \sigma\tilde r$
then
each blow-up limit at a point 
$(t^1,x^1)\in Q_{\tilde r}\cap \{ u=0\} \cap \{ \nabla u=0\}$ is time-independent.
\end{proposition}
\proof
Let us consider $(t^1,x^1)\in \{ u=0\} \cap \{ \nabla u=0\}$.
As the statement of the Proposition is by Theorem \ref{global-solution} true
when $(t^1,x^1)$ is a branch point, we may assume that $u\ge 0$
in some neighborhood of $(t^1,x^1)$.
From Lemma \ref{uniform-closeness} (with $\delta := \inf_{Q_1} \min(\lambda_+,\lambda_-)/(96n)$) and Proposition \ref{directional-monoton}
we know that $u$ is non-decreasing, say in the direction $e$
for every $e$ close to $x_n$ in $Q_{\tilde r}$
and that $|\partial_t u|\le \inf_{Q_{\tilde r}} \min(\lambda_+,\lambda_-)/2$
in $Q_{\tilde r}$.
\\
From \cite[Theorem 4.1]{siam} we infer now that
each blow-up limit $z$ at $(t^1,x^1)$ is a non-negative backward self-similar solution.
Concerning those, it has been
shown in \cite[Lemma 6.3]{cps} and \cite[Theorem 8.1]{cps} that
either $z$ is a half-plane solution of the form
$z(t,x) = \lambda_+(t^1,x^1)/2 \max(x\cdot e,0)^2$ for some
$e\in \partial B_1$, or
$z(t,x)=- a_0 t + \sum_{i=1}^n a_i x_i^2$ with non-negative 
constants $a_i, 0\le i\le n$ and $a_0\le \lambda_+(t^1,x^1)/2$.
In the latter case the symmetry of $z$ 
contradicts the fact that $z$ is non-decreasing in every
direction $e$ as above. Consequently $\partial_t z=0$ in $\{ t<0\}$,
and Lemma \ref{forward} and Corollary \ref{zerolim}
imply that
$\partial_t u(t^1,x^1)=0$.
\qed
\begin{corollary}\label{conti}
Let $u$ be a solution of (\ref{obst}) in $Q_1$.
Then
there are positive constants $\tilde r$ and $\tilde \sigma$
(depending on $\inf_{Q_{1}} \min(\lambda_+,\lambda_-)$,
the Lipschitz norms of $\lambda_\pm,$ the supremum norm of $u$ and
the space dimension $n$)
such that
the following holds.
If $u(0)=0\> ,\> |\nabla u(0)|\leq \tilde \sigma\tilde r$,
$\pardist(0,\{ u>0\})\le \tilde \sigma\tilde r$ and $\pardist(0,\{ u<0\})\le \tilde \sigma\tilde r$
then
$\partial_t u$ is continuous
in $Q_{\tilde r}$.
\end{corollary}
\proof
The corollary follows immediately from Lemma \ref{cinf}, Proposition \ref{time-indep}
and Corollary \ref{zerolim}.\qed
\section{Directional Monotonicity II}
It is now possible to extend the directional monotonicity result
of Section \ref{dirmon} to a directional monotonicity result
with respect to time-space variables.
\begin{proposition}\label{directional-monoton-par}
Let $0<\lambda_{\rm min}\le \inf_{Q_1(0)}
\min(1,\lambda_+,\lambda_-)$, $h$ as in (\ref{one-dim-sol}), let
$\varepsilon\in (0,1)$
and let $\tilde r$ and $\tilde \sigma$ be the constants of Corollary \ref{conti}. Then each solution $u$ of (\ref{obst}) in
$Q_1(0)$ such that
$$\dist_{W^{1,\infty}(Q_1(0))}(u,h)
\le \delta := {\lambda_{\rm min}\varepsilon\over {48 n}}\tilde r^2 \tilde \sigma^2$$ and
$$\sup_{Q_1(0)}\max(|\nabla \lambda_+|,|\partial_t \lambda_+|,|\nabla \lambda_-|,|\partial_t \lambda_-|)\le
\delta$$ satisfies $\varepsilon^{-1}\alpha \partial_t u \> + \> \varepsilon^{-1}
\partial_e u - |u| \ge 0$ in $Q_{1/2}(0)$ for every $\alpha \in [-1,1]$ and every $e\in \partial
B_1(0)$ such that $e_1\ge \varepsilon$;
here $e_1$ denotes the first component of the vector $e$.
\end{proposition}
\proof 
First note that $Q_1\cap\{ u=0\}$
is by the assumptions contained in the strip
$|x_1| < \tilde \sigma \tilde r/2$, implying by
Corollary \ref{conti} and Lemma \ref{cinf}
that $\partial_t u$ is continuous in $Q_1$. 
We know that $\varepsilon^{-1}\alpha\partial_t h \> +\> \varepsilon^{-1} \partial_e h - |h|\ge 0$.
It follows that
\begin{equation}\label{approx_p}
\varepsilon^{-1} \alpha\partial_t u \> +\> \varepsilon^{-1} \partial_e u - |u|\ge -3\delta \varepsilon^{-1}
\end{equation}
provided that $\dist_{W^{1,\infty}(Q_1(0))}(u,h) \le\delta$. Suppose now
towards a contradiction that the statement is not true. Then there
exist $\lambda_+,\lambda_-\in (\lambda_{\rm min},+\infty), (t^*,x^*) \in
Q_{1/2}(0), \alpha^*,e^*, $ and a solution $u$ of (\ref{obst}) in $Q_1(0)$
such that $\dist_{W^{1,\infty}(Q_1(0))}(u,h) \le\delta$,
$$\sup_{Q_1(0)}\max(|\nabla \lambda_+|,|\partial_t \lambda_+|,|\nabla \lambda_-|,|\partial_t \lambda_-|)\le
\delta,$$ $|\alpha^*|\le 1,e^*_1 \ge \varepsilon$ and
$\varepsilon^{-1}\alpha^* \partial_t u(t^*,x^*) \> +\>\varepsilon^{-1}\partial_{e^*} u(t^*,x^*) - |u(t^*,x^*)|<0.$ For the
positive constant $c$ to be defined later the functions $v :=
\varepsilon^{-1}\alpha^* \partial_t u \> +\>\varepsilon^{-1}\partial_{e^*} u - |u|$ and $w:=
\varepsilon^{-1}\alpha^* \partial_t u \> +\>\varepsilon^{-1}\partial_{e^*} u - |u|+c |x-x^*|^2 - c (t-t^*)$ satisfy then 
by the definition of $\delta$ the
following: in the set $D := Q_1(0)\cap \{ v<0\}\cap\{ t<t^*\}$,
\[ \Delta w \> -\> \partial_t w\le 2nc + c -{\lambda_+} \chi_{\{ u>0\}}-{\lambda_-} \chi_{\{ u<0\}}\]\[
+\varepsilon^{-1} (\lambda_+ +\lambda_-) \nu_x\cdot e^* {\mathcal
H}^{n}\lfloor (\{ u=0\}\cap \{ \nabla u\ne 0\}) \]\[
+\varepsilon^{-1} (\lambda_+ +\lambda_-) \nu_t \alpha^*  {\mathcal
H}^{n}\lfloor (\{ u=0\}\cap \{ \nabla u\ne 0\}) \]\[+\varepsilon^{-1}
( \chi_{\{ u>0\}}(\alpha^* \partial_t \> + \> \partial_{e^*}) \lambda_+ - \chi_{\{
u<0\}}(\alpha^* \partial_t \> + \>\partial_{e^*}) \lambda_-)
\]
where $\nu = {(\partial_t u,\nabla u) \over {|(\partial_t u,\nabla u)|}}$. As
\[\nu\cdot (\alpha^*,e^*)\le 0 \> \hbox{ on }\> \{ u=0\}\cap \{ v<0\} = \{ u=0\}\cap \{
\varepsilon^{-1}\alpha^* \partial_t u \> +\>\varepsilon^{-1}\partial_{e^*} u <0\}\; ,\] we obtain by the definition of $\delta$ that
$w$ is supercaloric in $D$ provided that $c$ has been chosen
accordingly, say $c:= \lambda_{\rm min}/(4n)$. It follows that the
negative infimum of $w$ is attained on
\[ \partial_{\rm par} D \subset (\partial_{\rm par} Q_1(0)\cap \{ t<t^*\})\cup
\left( Q_1(0)\cap \partial \{ v<0\}\right)
\; .\]
Consequently it is attained on $\{ t<t^*\}\cap \partial_{\rm par} Q_1(0)$,
say at the point $(\bar t,\bar x) \in \{ t<t^*\}\cap \partial_{\rm par} Q_1(0)$.
Since $\pardist((\bar t,\bar x),(t^*,x^*))\ge 1/2$, we obtain that
\[ \varepsilon^{-1}\alpha^* \partial_t u(\bar t,\bar x) \> +\>\varepsilon^{-1}\partial_{e^*} u(\bar t,\bar x) - |u(\bar t,\bar x)|\]\[ =
v(\bar t,\bar x) = w(\bar t,\bar x) - c |x^* - \bar x|^2 + c (\bar t-t^*) 
< -c/4=-\lambda_{\rm min}/(16n)\; .\] But this contradicts
(\ref{approx_p}) in view of $\delta = {\lambda_{\rm
min}\varepsilon\over {48 n}}\tilde r^2\tilde \sigma^2$.\qed
\section{Proof of the main theorem}
The theorem is proven in several simple steps, using mainly
Proposition \ref{directional-monoton-par}, and Lemma
\ref{uniform-closeness}. Note that the proof can be simplified substantially
in the case that we are dealing not with a whole class of solutions
but a single solution.
\newline
{\bf Part I:} In this first part we prove uniform Lipschitz regularity
and continuous differentiability with respect to the space variables.
\newline
{\bf Step 1 (Directional monotonicity):} Given $\varepsilon >0$, there are
$\sigma_\varepsilon>0$ and $r_\varepsilon >0$ (depending only on the parameters of the
statement) such that $2\alpha \varepsilon^{-1}r^2_\varepsilon\partial_t u+2\varepsilon^{-1}r_\varepsilon\partial_e u
-|u|\geq 0$ in $Q_{r_\varepsilon/2}(y)$ for every $\alpha \in [-1,1]$. The inequality holds for every $(s,y)\in
Q_{1/2}(0)$ satisfying
$u(s,y)=0$,$\>|\nabla u(s,y)|\leq \sigma_\varepsilon r_\varepsilon$,$\>\pardist((s,y),\{ u>0\})\le \sigma_\varepsilon r_\varepsilon$ and
$\pardist((s,y),\{ u<0\})\le \sigma_\varepsilon r_\varepsilon$,
for some unit vector
$\nu_\varepsilon(s,y)$ and for every $e \in \partial B_1$ satisfying $e\cdot \nu_\varepsilon(s,y)
\ge {\varepsilon \over 2}$. In particular, for $\varepsilon=1$, the solution $u$ is 
by condition (\ref{cond}) with $\sigma=\sigma_1 r_1$
non-decreasing in
$Q_{r_1/2}(0)$ in direction $(r_1,e)$ for every $e \in \partial B_1(0)$ such that $e \cdot
\nu_\varepsilon(0) \geq {1\over 2}$.
\newline
{\it Proof:} By Lemma \ref{uniform-closeness} there are
$\sigma_\varepsilon>0$ and $r_\varepsilon>0$ as above such that the scaled function
$u_{r_\varepsilon}(t,x) = u(s+r_\varepsilon^2 t,y+r_\varepsilon x)/r_\varepsilon^2$ is
$\delta:=\varepsilon \frac{\lambda_{\rm min}}{64n}\tilde r^2\tilde \sigma^2$-close in
$C^1(Q_1(0))$ to a rotated version $\tilde h$ of $h$ in $Q_1$. 
Let $\nu_\varepsilon(s,y)$ be the accordingly rotated version of the unit
vector $e_1$.
Since
$u_{r_\varepsilon}$ solves (\ref{obst}) with respect to
$\lambda_+(r_\varepsilon^2 \cdot + s, r_\varepsilon \cdot +y)$ and $\lambda_-(r_\varepsilon^2 \cdot + s,r_\varepsilon
\cdot +y)$, and since  $\max(|\nabla (\lambda_+(r_\varepsilon^2 \cdot + s,r_\varepsilon \cdot
+ y))|,|\nabla (\lambda_-(r_\varepsilon^2 \cdot + s,r_\varepsilon \cdot+y))|,|\partial_t (\lambda_+(r_\varepsilon^2 \cdot + s,r_\varepsilon \cdot
+ y))|,|\partial_t (\lambda_-(r_\varepsilon^2 \cdot + s,r_\varepsilon \cdot
+ y))|)\le C_1
r_\varepsilon $, we may choose $r_\varepsilon < \delta/C_1$ in order
to apply
 Proposition \ref{directional-monoton-par} to
$u_{r_\varepsilon}$ in $Q_1$ and to conclude that
$2\alpha \varepsilon^{-1}\partial_t u_{r_\varepsilon}+2\varepsilon^{-1}\partial_e u_{r_\varepsilon} - |u_{r_\varepsilon}|\ge 0$ in $Q_{1/2}(0)$ for
every $\alpha \in [-1,1]$ and every $e\in \partial B_1(0)$ such that $e\cdot \nu_\varepsilon(s,y) \ge
\varepsilon/2$. Scaling back we obtain the statement of Step 1.
\newline
{\bf Step 2 (Lipschitz continuity):} $\partial \{u>0\}\cap
Q_{r_1/2}(0)$ and $\partial \{u<0\}\cap Q_{r_1/2}(0)$ are Lipschitz
graphs in the direction of $(0,\nu_\varepsilon(0))$ with spatial Lipschitz norms less than
$1$ and temporal Lipschitz norms less than $r_1^{-1}$. Moreover, for each $\varepsilon\in (0,1)$ and $(s,y)\in
\{ u=0\} \cap Q_{1/2}$ satisfying $|\nabla u(s,y)|\leq \sigma_\varepsilon r_\varepsilon$,
$\pardist((s,y),\{ u>0\})\le \sigma_\varepsilon r_\varepsilon$ and
$\pardist((s,y),\{ u<0\})\le \sigma_\varepsilon r_\varepsilon$,
the free boundaries $\partial
\{ u>0\} \cap Q_{r_{\varepsilon /2}}(s,y)$ and $\partial \{ u<0\} \cap
Q_{r_{\varepsilon /2}}(s,y)$ are Lipschitz graphs (in the direction of
$\nu_\varepsilon(s,y)$) with spatial Lipschitz norms not greater than $\varepsilon$.
\newline
{\it Proof:} This follows from the monotonicity obtained in Step 1.
\newline
{\bf Step 3 (Existence of a spatial tangent plane at points
$(s,y) \in \partial\{ u>0\}\cap \partial \{ u<0\}\cap Q_{1/2}(0)$ satisfying $|\nabla u(s,y)|=0$):}
The Lipschitz graphs of Step 2 are both
differentiable with respect to the space variables
at the point $(s,y)$, and the two spatial tangent planes
at $(s,y)$ coincide.
\newline
{\it Proof:} This follows from Step 2 by letting $\varepsilon$ tend
to zero.
\newline
{\bf Step 4 (One-phase points are regular):} If $(s,y) \in Q_{r_1/2}(0)$
is a free boundary point and the solution $u$ is non-negative or
non-positive in $Q_\delta(s,y)$, then the free boundary is the graph
of a $C^{1,\alpha}$-function in $Q_{c_1\delta}(s,y)$, where $c_1$ and
the $C^{1,\alpha}$-norm depend only on the parameters in the
statement. Consequently, in $Q_{r_1/2}(0)$, there exist no singular
one-phase free boundary points.
\newline
{\it Proof:} By Step 2, the sets
$\{ u>0\} \cap Q_{r_1/2}(0)$ and $
\{ u<0\} \cap Q_{r_1/2}(0)$ are sub/supergraphs of
Lipschitz continuous functions.
Therefore $\{ u=0\}\cap Q_\delta(s,y)$
satisfies the thickness condition required for
\cite[Theorem 15.1]{cps}
and the statement follows.
\newline
{\bf Step 5 (Existence of space normals in $Q_{r_1/2}(0)$):}
$\partial
\{ u>0\} \cap Q_{r_1/2}(0)$ and $\partial
\{ u<0\} \cap Q_{r_1/2}(0)$ are graphs
of Lipschitz continuous functions which are differentiable
with respect to the space variables. 
\newline
{\it Proof:} Let $(s,y) \in Q_{r_1/2}(0)$ be a free boundary point.
We have to prove existence of a tangent plane at $(s,y)$.\newline
First, if $(s,y)$ is a one-phase point,
i.e. if the solution $u$ is non-negative or non-positive in
$Q_\delta(s,y)$, then the statement holds at $(s,y)$ by the result of Step 4.
Second, if $|\nabla u(s,y)|\ne 0$, the statement holds by Lemma \ref{cinf}.
Last, if $|\nabla u(s,y)|=0$ and $(s,y)$ is the limit point
of both phases $\{ u>0\}$ and $\{ u<0\}$, then
Step 3 applies.
\newline
{\bf Step 6 (Equicontinuity of the space normals):} It remains to prove
that the space normals are equicontinuous on $Q_{r_1/2}(0)\cap \partial \{
u>0\}$ and on $Q_{r_1/2}(0)\cap \partial \{ u<0\}$ for $u$ in the
class of solutions specified in the statement of the main theorem.
\newline\noindent
{\it Proof:} By Step 2 we know already that the spatial Lipschitz norms of
$\partial \{ u>0\} \cap Q_{r_1/2}(0)$ and $\partial \{ u<0\} \cap
Q_{r_1/2}(0)$ are less than $1$. We prove that the space normals are
equicontinuous on $Q_{r_1/2}(0)\cap\partial \{ u>0\}$.
\newline
We may assume that $\nu(0)$ points in the direction of the $x_1$-axis and
that $x_1=f(t,x_2,\dots,x_n)$ is the representation of $\partial
\{u>0\}\cap Q_{r_1/2}(0)$. Besides we have $|\nabla f(t,x')|<1$ for
$(t,x)=(t,x_1,x') \in
\partial \{u>0\}\cap Q_{r_1/2}(0)$.
We claim that for $\varepsilon >0$ there is $\delta_\varepsilon >0$
depending only on the parameters in the statement such that for any
pair of free boundary points $(s^1,y^1),(s^2,y^2) \in
\partial \{ u>0\} \cap Q_{r_1/2}(0)$,
\begin{equation}\label{continuity-of-normal}
\pardist((s^1,y^1),(s^2,y^2)) \leq \delta_\varepsilon \quad \Rightarrow \quad |\nu(s^1,y^1)
-\nu(s^2,y^2)|\leq 2\varepsilon .
\end{equation}
In what follows let $\rho_\varepsilon := \sigma_\varepsilon r_\varepsilon/2 \le r_1/2$.
\newline
Suppose first that $u$ is non-negative in $Q_{\rho_\varepsilon}(s^1,y^1)$. Here we
may as in Step 4 apply \cite[Theorem 15.1]{cps} to the scaled
function $w(t,x) := u(s^1+\rho_\varepsilon^2 t, y^1+\rho_\varepsilon x)/{\rho^2_\varepsilon}$; since
the $C^{1,\alpha}$-norm of the free boundary normal of $w$ is on $Q_{c_2}\cap \partial\{ w>0\}$
bounded by a constant $C_3$, where $c_2>0$ and $C_3<+\infty$ depend only on the parameters in the
statement, we may choose
$$\delta_\varepsilon :=  \min({\varepsilon^{1 \over \alpha}\over {C_3^{1 \over \alpha}}},c_2)\rho_\varepsilon$$
to obtain (\ref{continuity-of-normal}).
\newline
Next, suppose that 
$u$ changes its sign at $Q_{\rho_\varepsilon}(s^1,y^1)$.
If there is a point
$(s,y)\in
Q_{\rho_\varepsilon}(s^1,y^1)\cap \partial\{ u>0\}$ such that
$|\nabla u(s,y)|\leq\rho_\varepsilon$ then we are in the situation of Step 1. By Step 2 the
free boundary $\partial \{ u>0\} \cap Q_{r_{\varepsilon}/2}(s,y)$ is
Lipschitz with spatial Lipschitz norm not greater than ${\varepsilon}$. 
Hence (\ref{continuity-of-normal})
follows in this case with $\delta_\varepsilon := r_{\varepsilon}/2$.
\newline
Last, if $|\nabla u(s,y)|\geq \rho_\varepsilon$
for all points $(s,y)\in Q_{\rho_\varepsilon}(s^1,y^1)\cap \partial\{ u>0\}$,
we proceed as follows:
from the equation $u(t,f(t,x'),x')=0$ we infer that
$\nabla' u+\partial_1 u\> \nabla'f=0$ 
on $\partial \{u>0\}\cap
Q_{r_1/2}(0)$. Hence we obtain
$$
|\nabla' f(s^1,(y^1)')-\nabla' f(s^2,(y^2)')|
= \left\vert {\nabla' u(s^1,y^1)\over {\partial_1 u(s^1,y^1)}}-
{\nabla' u(s^2,y^2)\over {\partial_1 u(s^2,y^2)}}\right\vert$$
$$\le {|\nabla' u(s^2,y^2)-\nabla' u(s^1,y^1)|\over {|\partial_1 u(s^1,y^1)|}}$$
$$+ \> \left|{\nabla' u(s^2,y^2)\over {\partial_1 u(s^2,y^2)}}\right|
{|\partial_1 u(s^2,y^2)-\partial_1 u(s^1,y^1)|\over {|\partial_1 u(s^1,y^1)|}}$$
$$ \leq 4M\rho_\varepsilon^{-1}\pardist((s^1,y^1),(s^2,y^2))\; ,
$$
where $M=\Vert \nabla u \Vert_{C^{1/2,1}(Q_{1/2}(0))}$.
In particular we may choose
$$\delta_\varepsilon :=  {\varepsilon \over {4M}}\rho_\varepsilon$$
to arrive at (\ref{continuity-of-normal}).\\[.3cm]
{\bf Part II:} We conclude the proof of the main theorem by pointing out a counter-example to $C^1$-regularity.
\newline
\begin{figure}
\begin{center}
\input{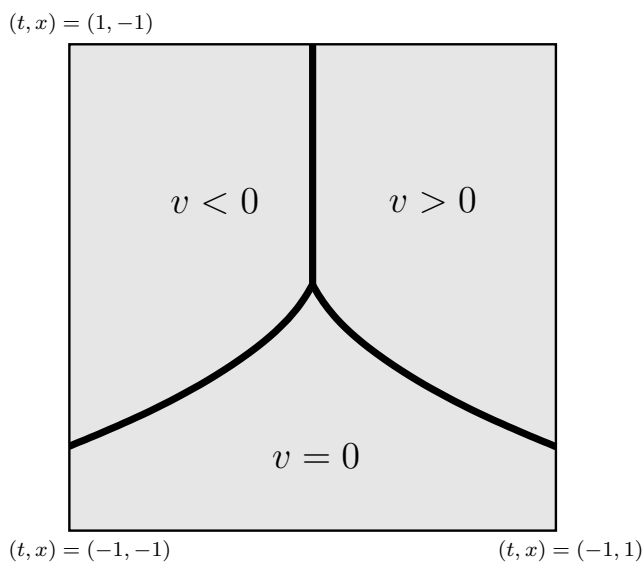}
\end{center}
\caption{A counter-example to $C^1$-regularity}\label{counter}
\end{figure}
Consider the one-phase counter-example $u:[-r^2,r^2]\times [0,r]\to [0,+\infty)$
from \cite[p. 376]{asu} satisfying the following:
$u(t,0)=0$ for $-r^2\le t \le r^2$,
and the free boundary touches the lateral boundary at the origin
in a non-tangential way.
Thus we may reflect $u$ to a solution 
$$v(t,x) :=\left\{ \begin{array}{l}
u(t,x),x\ge 0\\
-u(t,-x),x<0\end{array}\right.$$
and obtain that $v$ is a solution of our two-phase problem
(\ref{obst}) in $Q_r$ for $\lambda_+=\lambda_-=1$.
As the free boundary $\partial\{ v>0\}$ is only Lipschitz at the origin,
we conclude that differentiability with respect to the time variable is 
in general not true.\qed
\bibliographystyle{plain}
\bibliography{suw_par_final.bib}
\end{document}